\documentclass[12pt]{article}

\usepackage[utf8]{inputenc}
\usepackage[T1]{fontenc}
\usepackage{graphicx}
\usepackage{longtable}
\usepackage{wrapfig}
\usepackage{rotating}
\usepackage[normalem]{ulem}
\usepackage{amsmath}
\usepackage{amssymb}
\usepackage{capt-of}
\usepackage{hyperref}
\usepackage{amsthm,amssymb,amsmath}
\usepackage{hyperref}
\usepackage{graphicx}
\usepackage{enumitem}
\usepackage{tabu}
\theoremstyle{definition}
\newtheorem{definition}{Definition}
\newtheorem{mexample}{Example}
\theoremstyle{remark}
\newtheorem{remark}{Remark}
\theoremstyle{plain}
\newtheorem{lemma}{Lemma}

\newtheorem{theorem}[lemma]{Theorem}
\newtheorem{corollary}[lemma]{Corollary}
\newtheorem{conjecture}{Conjecture}
\newtheorem{question}{Question}
\newcommand{\modd}[2]{#1\ \textup{(mod }#2\textup{)}}
\newcommand{\set}[1]{\left\{{#1}\right\}}

\newcommand\setsuchas[2]{\left\{{#1}:{#2}\right\}}
\newcommand{\NN}{{\mathbb{N}}}
\newcommand{\NNp}{{\mathbb{P}}}
\newcommand{\ZZ}{{\mathbb{Z}}}
\newcommand{\CC}{{\mathbb{C}}}

\newcommand{\order}{{\mathrm{o}}}
\newcommand{\divisors}[1]{{D(#1)}}
\newcommand{\type}{\tau}
\newcommand{\Arank}[1]{r_{#1}}
\newcommand{\gdecomp}[1]{\Pi^{(#1)}}
\newcommand{\gpartition}[1]{\mathcal{S}^{(#1)}}
\newcommand{\gconvolution}[1]{\mathcal{G}^{(#1)}}
\newcommand{\lcm}{\operatorname{lcm}}
\newcommand{\val}{\operatorname{val}}
\newcommand{\divides}[2]{#1 \mid #2}
\newcommand{\dividesnot}[2]{#1 \nmid #2}
\newcommand{\height}[2]{\mathbf{h}_{#1}(#2)}
\newcommand{\ssift}[1]{\mathbf{s}({#1})}
\newcommand{\seqnum}[1]{\href{https://oeis.org/#1}{\textrm \underline{#1}}}
\newcommand{\dens}{\mathbf{d}}
\newcommand{\norm}[1]{\lVert #1 \rVert}
\author{Jan Snellman \\ Department of Mathematics, Link{\"o}ping University \\ 58183 Link{\"o}ping, Sweden \\ \texttt{jan.snellman@liu.se}}
\date{2025-2026}
\title{Greedy Regular Convolutions}
\hypersetup{
 pdfauthor={},
 pdftitle={Greedy Regular Convolutions},
 pdfkeywords={arithmetic functions, regular convolutions, sieve methods},
 pdfsubject={},
 pdfcreator={},
 pdflang={English}}
\begin{document}

\maketitle

\begin{abstract}
Narkiewicz classified  \emph{regular convolutions} on arithmetical functions
by means of families \(\{\pi_{p} : p \text{ prime}\}\) of decompositions of \(\NN\) into arithmetic
progressions. Most regular convolutions, such as the \emph{Dirichlet convolution}
and the \emph{unitary convolution} are \emph{homogeneous}:  all \(\pi_{p}\) are equal.
The unitary convolution is furthermore \emph{bounded}: there is a common finite bound for
the lengths of the blocks  of \(\pi_{p}\).

In this article we define, for each positive integer \(d\),  the \emph{greedy convolution} of length \(d\).
It is a regular, homogeneous, bounded convolution where the maximal length of a block of \(\pi_{p}\) is \(d\).
The case \(d=1\) yields the unitary convolution and the case \(d=2\) produces the "ternary convolution".
These two convolutions are the only regular, homogeneous convolutions where
all blocks of \(\pi_{p}\) have the same length.

Thus, for the greedy convolutions with block-length \(d > 2\), the blocks of \(\pi_{p}\) have unequal length.
We describe the case \(d=3\) in detail: the lengths are 3 and 1.
The set of \emph{primitive exponents}, the smallest positive element of each block in the decomposition \(\pi_{p}\),
can be generated by a simple recursive procedure that we name \emph{selective sifting}.

The structure of the primitive exponents for greedy convolutions of higher length is surprisingly
intricate, given the straight-forward definition. We give a conjecture about the case \(d=4\).
\end{abstract}
\section{Introduction}
\label{sec:orged176b9}
\subsection{Arithmetical functions and regular convolutions; Narkiewicz classification}
\label{sec:org5cc9cc8}
We denote the non-negative integers with \(\NN\) and the positive integers with \(\NNp\).
The set of \emph{arithmetical functions}, i.e., complex-valued functions defined on \(\NNp\),
is denoted by \(\Gamma\).
Some distinguished elements of \(\Gamma\) are
the \emph{zeta function}
\begin{equation}
\label{eq:orgb700fb9}
  \begin{split}
     \zeta : \NNp & \to \CC \\ 
     n & \mapsto 1,    
  \end{split}
\end{equation}
the \emph{zero function}
\(\mathbf{0}\). For each positive integer \(j\) the
\emph{indicator function}
\begin{equation}
\label{eq:org0e7e65f}
e_{j}(n) =
\begin{cases}
1, & n = j; \\
0, & n \neq j.
\end{cases}
\end{equation}
is an arithmetical function. For each prime \(p\), the \(p\)-valuation
\begin{equation}
\label{eq:org109925e}
\val_{p}(n) = \max \setsuchas{a}{\divides{p^{a}}{n}}
\end{equation}
is an arithmetical function.

Arithmetical functions can be multiplied using various convolution products.
The most widely used such product is the \emph{Dirichlet convolution} \cite{vaidy,Afunc},
\begin{equation}
\label{eq:orgeb87452}
(f * g) (n) = \sum_{\substack{(a,b) \in \NNp \times \NNp \\ ab = n}} f(a)g(b) ,
\end{equation}
and another common variant is the \emph{unitary convolution}
\cite{Schinzel:Property,UniDiv,UniProd,Snellman:UniDivTop,Snellman:UniTrunc},
\begin{equation}
\label{eq:org182c43e}
(f * g) (n) = \sum_{\substack{(a,b) \in \NNp \times \NNp \\ab = n \\ \gcd(a,b)=1}} f(a)g(b).
\end{equation}
The \emph{ternary convolution}, introduced by Gavell in her Bachelor's thesis \cite{Gavell:Regular}, is defined
by the following daunting formula (a simpler formulation will be given later):
\begin{equation}
\label{eq:org9b5b66e}
  (f * g) (n) =
  \sum_{\substack{
    (a,b) \in \NNp \times \NNp \\
    ab = n \\
    \left(p \text{ prime and } \divides{p}{a} \text{ and } \divides{p}{b}\right) \implies \left( \val_p(a) \text{ even and} \val_p(b) \text{ even} \right)}}
  f(a)g(b)
\end{equation}

Narkiewicz \cite{Nark:conv} classified \emph{regular convolutions},
a very general family of convolutions that includes the Dirichlet, unitary, and ternary convolutions.
These regular convolutions yields a multiplication with some natural and desired properties;
commutativity and associativity, among others.
They are constructed by choosing, for each prime \(p\), a
"decomposition into arithmetic progressions of \(\NN\) containing zero, having
no positive integer in common"  \(\pi_{p}\).
Each \emph{block} in \(\pi_{p}\) contains the element 0;
so forming \(\tau_{p} = \setsuchas{s \setminus \set{0}}{s \in \pi_p}\) gives a set-partition of \(\NNp\).
This also partitions
\(\setsuchas{p^{a}}{a \in \NNp}\), and this partition is used to define the multiplication \((f*g)(n)\).
We will recall this construction, and Narkiewicz classification, in the next section.

All the regular convolutions introduced so far uses the same partition of \(\NNp\) for
all primes \(p\); we call such regular convolutions \emph{homogeneous}.
Furthermore, while the Dirichlet convolution has a single, infinite part in this partition, the
parts in the partition associated with the unitary convolution are all singletons, hence finite. We
call a regular convolution \emph{bounded}  if there is a common finite bound for the sizes of the parts of \(\tau_{p}\).
The ternary convolution is bounded, as well.

When \(A = \set{0, k, 2k, \dots, dk}\) is a finite block, we say that \(A\) has \emph{block-length} \(d\);
this is of course the cardinality of the  corresponding part
\(A \setminus \set{0} = \set{k, 2k, \dots, dk} \in \pi_p\). 

In Table \ref{tab:orge53dce4} we show the "blocks" of \(\pi_{p}\) and the "parts" of \(\tau_{p}\),
for the unitary, the ternary, and the greedy convolution of block-length \(p\).
These convolutions are homogeneous, so  \(p\) is arbitrary.
The greedy convolutions will be defined shortly.

\begin{table}[htbp]
\caption{\label{tab:orge53dce4}The first parts and blocks of the unitary, ternary, and the greedy convolution of block-length 3.}
\centering
\begin{tabular}{llllllll}
convolution & p/b &  &  &  &  &  & \\
\hline
unitary & parts & 1 & 2 & 3 & 4 & 5 & 6\\
unitary & blocks & 0,1 & 0,2 & 0,3 & 0,4 & 0,5 & 0,6\\
ternary & parts & 1,2 & 3,6 & 4,8 & 5,10 & 7,14 & 9,18\\
ternary & blocks & 0,1,2 & 0,3,6 & 0,4,8 & 0,5,10 & 0,7,14 & 0,9,18\\
greedy 3 & parts & 1,2,3 & 4,8,12 & 5,10,15 & 6 & 7,14,21 & 9,18,27\\
greedy 3 & blocks & 0,1,2,3 & 0,4,8,12 & 0,5,10,15 & 0,6 & 0,7,14,21 & 0,9,18,27\\
\end{tabular}
\end{table}

One can ask: what are the regular, homogeneous, bounded convolutions where all block-lengths
are equal, and finite? The unitary convolution have blocks with block-length 1.
Gavell \cite{Gavell:Regular} proved that there is a unique regular convolution having all block-lengths 2,
namely the ternary convolution,
but that for \(d > 2\)
there is no regular convolution having all block-lengths equal to \(d\).

What if we settle for blocks of length \(\le d\), filling
the blocks as much as possible? The \emph{greedy convolutions} introduced in this work is one answer to
this question. Here, we successively add the integers of \(\NNp\), which we must partition into
arithmetic progressions of length \(\le d\), into the \textbf{first available block with free space}.
As an example, in Figure \ref{fig:orga3a2aae} the integers 1 to 40 has been added; some blocks are final,
such as all of block-length 3, as well as the block \(\set{0,6}\),
some are yet to be finished. The integer 12 could extend two arithmetic progressions: \(\set{0,4,8}\) and
\(\set{0,6}\); the greedy rule places it in the first block, rather than in the second.

This procedure will yield the unitary convolution for \(d=1\) and the ternary convolution for
\(d=2\); in those cases, there is no choice involved.
We will study the case \(d=3\) in some detail in this article; the blocks are as shown in
Table \ref{tab:orge53dce4}.

Of special interest are the \emph{primitive} elements (or exponents), by which we mean the smallest element
in each part (smallest positive element in each block). The multiplication of arithmetical functions
using a homogeneous regular convolution is determined by the multiplication table of the indicator functions
\(e_{{p^{b}}}\), for \(p^{b}\) a prime power with \(b\) primitive. 

For the greedy convolutions of small block-length, we describe the set of primitives
using a recursive procedure that we have named
\emph{selective sifting}. This procedure produces, given a set \(S\) of integers \(\ge 2\), the subset
\(\ssift{\set{S}} \subset \NNp\) of those \(n\) which are not "selectively sifted out" by removing \(n\)
if \(n=ab\) with \(a \in S\), \(b \in \ssift{\set{S}}\). The set of primitives
for the ternary convolution consists of all positive integers with an even number of twos in their
prime factorization; this is \(\ssift{\set{2}}\). The primitives for the greedy convolution with \(d=3\)
are \(\ssift{\set{2,3,6}}\).

For greater block-lengths, the structure of the primitives becomes very complicated;
we offer a vague conjecture (involving selective sifting) for block-length four.
\subsection{Outline of the article}
\label{sec:orgc5989c6}
We give a brief review of the vector space \(\Gamma = \CC^{\NNp}\) of arithmetical functions in section
\ref{sec:org5b591e8}, recalling the usual topology,
which turns \(\Gamma\) into an \emph{ultrametric} space.
Next, we consider convolutions (multiplications) on \(\Gamma\).
The class of \(A\)-convolutions \cite{BurnettGeneralAfunctions} is
mentioned, but we focus on the class of \emph{regular convolutions} introduced and classified
by Narkiewicz \cite{Nark:conv}. In particular, we discuss \emph{homogeneous} and \emph{bounded}
regular convolutions, which are defined by a single partition of \(\NNp\) into
finite arithmetic progressions. We recall the notions of \emph{primitive powers}, \emph{rank}, and \emph{height}.
For convenience, we introduce the notions of \emph{primitive elements}.

The ternary convolution is discussed in section \ref{sec:orgeaefe31}. In particular,
we describe the structure of the set of its primitive elements, and  calculate
the \emph{natural density} of that set as a subset of \(\NNp\).

Section \ref{sec:org603cab4} introduces the class of \emph{greedy convolutions}. We give examples
and show some simple properties. We also present Gavell's \cite{Gavell:Regular} proof
of the fact that the unitary and the ternary convolution are the only regular, homogeneous,
bounded convolutions with equal block-length.

The case \(d=3\) of greedy convolutions is studied in section \ref{sec:org1f75020}.
The primitive elements are described, and the natural density of that set is calculated.

In section \ref{sec:org6b6b851} we introduce the process of \emph{selective sifting}, a sieve-like
recursive method that provides a compact description of the set of primitives for the
greedy convolution of length 2 and 3, and which seems to be pertinent for greedy convolutions
of higher block-length, as well. The procedure is curious and amusing and might be interesting
in other contexts. We give some simple result, and pose some questions for further investigations.

The greedy convolution of block-length 4 is briefly discussed in
section \ref{sec:orgfd29dc1}. We present the vaguest of conjectures for
the primitive elements, backed by some experimental data.

Finally, in section \ref{sec:org44d0849} we collect a number of
open problems about the primitive elements of greedy convolutions of block-length \(d > 3\).
\section{Convolution products on arithmetical functions}
\label{sec:org5b591e8}
\subsection{The vector space \(\Gamma\) and its topology}
\label{sec:org0cede2a}
The set of all \emph{arithmetical functions}, complex-valued functions defined on the positive integers, is
denoted \(\Gamma = \CC^{\NNp}\). It is a \(\CC\)-vector space.
It is also a topological space when endowed with the
product topology. Equivalently, we can define the \emph{order} \(\order(f)\) of \(f \in \Gamma\) as
\begin{equation}
\label{eq:orgbf47533}
\order(f) = \sup \setsuchas{j \in \NNp}{i < j \implies f(i) = 0},
\end{equation}
and then the \emph{norm} of \(\norm{f}\) as \(1/\order(f)\); the \emph{distance} between two elements of \(\Gamma\)
is then
\[d(f,g) =\norm{f-g}.\]
This metric is an \emph{ultrametric};
\[d(f,h) \le \max \left( d(f,g), d(g,h) \right),\]
and the topology induced by this ultrametric is the same as the product topology \cite{vring,dvar}. 

For any \(f \in \Gamma\),
\begin{equation}
\label{eq:org5a39e0b}
f = \sum_{j \in \NNp} f(j)e_{j},
\end{equation}
and this sum converges absolutely, irrespective of reordering of terms.

\begin{remark}
In this context, the ring \(\CC\) is
usually given the discrete topology.
This means that  the sequence of constant functions \(2^{-k} \zeta\) does not
converge to the zero function; rather, each function in the sequence have norm 1.
\end{remark}

An arithmetical function \(f \in \Gamma\) is \emph{multiplicative} if \(f(mn)=f(m)f(n)\) whenever \(\gcd(m,n)=1\),
and \emph{totally multiplicative} if \(f(mn)=f(m)f(n)\) for all \(m,n\).
The zeta function is totally multiplicative.
\subsection{\(A\)-convolutions on \(\Gamma\)}
\label{sec:orgc10c7e3}
There are numerous well-studied \emph{convolution products} on \(\Gamma\) that turns it into a \(\CC\)-algebra.
Such a convolution is usually required to be associative, \(\CC\)-bilinear, and continuous as a
function \(\Gamma \times \Gamma \to \Gamma\). 

We have already defined the Dirichlet convolution (\ref{eq:orgeb87452}), the unitary convolution
(\ref{eq:org182c43e}), and the ternary convolution (\ref{eq:org9b5b66e}). They all fit into a common
framework that covers many, but certainly not all, convolutions, and which we now describe. Our
exposition follows \cite{BurnettGeneralAfunctions}.

\begin{definition}
Let \(\divisors{n}\) denote the set of positive divisors of the positive integer \(n\).
Suppose that
\begin{equation}
\label{eq:orgb56f580}
\begin{split}
  A: \NNp & \to 2^{\NNp}\\
  A(n) & \subset \divisors{n}
\end{split}
\end{equation}
is a function which gives, for each positive integer \(n\), a subset \(A(n)\) of the
positive divisors of \(n\).
Then, given \(f,g \in \Gamma\), we define the "\(A\)-convolution"
\begin{equation}
\label{eq:org1d03035}
(f*_{A}g)(n)= \sum_{d \in A(n)}f(d)g(\frac{n}{d}).
\end{equation}
\label{orgc6bd955}
\end{definition}

\begin{definition}
Let \(A\) be as above. We say that \(A\) (or its associated convolution \(*_{A}\)) is
\begin{enumerate}[label={(\alph*)}]
\item \emph{simple} if \(1 \in A(n)\) for all \(n\),
\item \emph{reflexive} if \(n \in A(n)\) for all \(n\),
\item \emph{symmetric} if \(d \in A(n)\) implies that \(\frac{n}{d} \in A(n)\),
\item \emph{transitive} if \(d \in A(n)\) implies that  \(A(d) \subseteq A(n)\),
\item \emph{multiplicative} if \(A(mn) = A(m)A(n)\) whenever \(\gcd(m,n)=1\); here
\[A(m)A(n) = \setsuchas{uv}{u \in A(m), v \in A(n)},\]
\item \emph{homogeneous} if \(A\) is multiplicative, and for every pair of primes \(p_{1},p_{2}\) and
every \(0 \le b \le a\), \(p_{1}^{b} \in A(p_{1}^{a})\) if and only if \(p_{2}^{b} \in A(p_{2}^{a})\).
\end{enumerate}
\label{org53816d8}
\end{definition}

The important class of \emph{regular} convolutions is defined as follows.
\begin{definition}
The convolution \(*_{A}\) is \emph{regular} if the following 3 conditions hold:
\begin{enumerate}[label={(\roman*)}]
\item the ring \((\Gamma,+,*_A)\) is associative, commutative, and possesses a neutral element for multiplication,
\item whenever \(f,g \in \Gamma\) are multiplicative arithmetical functions, then so is \(f *_{A} g\),
\item the zeta function \(\zeta\) has a unique two-sided inverse \(\mu_{A}\) such that
\(\mu_{A}(p^{a}) \in \set{0,-1}\) for all prime powers \(p^{a}\).
\end{enumerate}
\label{org5fb700c}
\end{definition}

Narkiewicz \cite{Nark:conv} classified the regular convolutions.
\begin{theorem}[Narkiewicz]
The convolution \(*_{A}\) is regular if and only if it has the following two properties:
\begin{enumerate}[label={(\roman*)}]
\item \(A\) is simple, reflexive, symmetric, transitive, and multiplicative,
\item for each prime power \(p^{k}\), there is some positive divisor \(t\) of \(k\)
such that
\begin{equation}
\label{eq:org559f6b8}
 A(p^{k}) = \set{1, p^{t}, p^{2t}, \dots, p^{rt}=p^{k}}, 
\end{equation}
and furthermore
\[
  p^{t} \in A(p^{2t}), \, p^{2t} \in A(p^{3t}), \,  \cdots , \,  p^{(r-1)t} \in A(p^{k}).
  \]
\end{enumerate}
\label{org0aec854}
\end{theorem}

\begin{theorem}[Narkiewicz]
Let \(K\) denote the class of all decompositions of \(\NN\) into
arithmetic progressions (finite or not) containing zero, such that no two progressions
belonging to the same decomposition have a positive integer in common. Associate to
each prime \(p\) the decomposition \(\pi_{p} \in K\), and define
\begin{equation}
\label{eq:org40e1f9c}
  \prod p_{i}^{a_{i}} \in A \Bigl( \prod p_{i}^{b_{i}} \Bigr) \quad \iff \quad \forall i: a_{i} \le b_{i}
  \text{ and } a_{i}, b_{i} \text { belong to the same block of } \pi_{p_{i}}.
\end{equation}
Then \(A\) is regular, and any regular convolution is obtained this way.
\label{org0791e1e}
\end{theorem}

Narkiewicz defined the notions of \emph{primitiveness}, \emph{type}, and \emph{rank} for prime powers:
\begin{definition}
For a regular convolution \(A\), and a prime power \(p^{k}\), suppose that
\begin{equation}
\label{eq:org88b7a90}
A(p^{k}) = \set{1,p^{t}, p^{2t}, \dots, p^{rt}=p^{k}}.
\end{equation}
Then \(t = \type_{A}(p^{k})\)  is the \emph{type} of \(p^{k}\).
We say that \(p^{k}\) is \emph{primitive} if \(A(p^{k}) = \set{1,p^{k}}\).
If \(p^{k}\) is primitive, the \emph{rank}  of \(p^{k}\) is
\[
\Arank{}(p^{k}) = \Arank{A}(p^{k})=  \mathrm{sup} \setsuchas{r}{p^{k} \in A(p^{r})}.
\]
\label{org5d74b9d}
\end{definition}

\begin{remark}
Note that
\begin{displaymath}
e_{a} *_{A} e_{b} =
\begin{cases}
  e_{ab}, & a,b \in A(ab); \\
  \mathbf{0}, & \text{ otherwise. }
\end{cases}
\end{displaymath}
This---together with distributivity, \(\CC\)-linearity, and continuity---determines the regular convolution \(*_{A}\), in the sense that
\begin{equation}
\label{eq:org9aee4e2}
  \Bigl(\sum_{i} f(i)e_{i}  \Bigr) *_{A} \Bigl(\sum_{j} g(j)e_{j}  \Bigr) =
  \sum_{i} \sum_{j} f(i)g(j) (e_{i}*_{A} e_{j}).
\end{equation}
In fact, the regular convolution \(*_{A}\) is determined by
\begin{equation}
\label{eq:org3616d57}
e_{p^{a}} *_{A} e_{p^{b}} =
\begin{cases}
  e_{p^{a+b}}, & a,b,a+b \text{ in same block of } \pi_{p}; \\
  \mathbf{0}, & \text{ otherwise, }
\end{cases}
\end{equation}
since \(e_{m} *_{A} e_{n} = e_{mn}\) whenever \(\gcd(m,n)=1\). Furthermore, the functions \(e_{{p^{k}}}\) with \(p^{k}\) a primitive prime power
constitute an algebraically independent set, generating a dense subring of \((\Gamma,+,*_{A})\).
\label{orgb245134}
\end{remark}
\subsection{Homogeneous bounded regular convolutions}
\label{sec:org1d3d116}
An \(A\)-convolution is \emph{homogeneous} (in the sense of Burnett and Taylor \cite{BurnettGeneralAfunctions})
if for every pair of primes \(p_{1},p_{2}\), and every pair of non-negative integers \(a,b\) with \(0 \le b \le a\),
\[p_{1}^{b} \in A(p_{1}^{a}) \iff p_{2}^{b} \in A(p_{2}^{a}).\]
If \(A\) is homogeneous and regular, it is determined by a single partition
\(S\) of \(\NNp\), having the property that \(s \cup \set{0}\) is an arithmetic progression for
all \(s \in S\). We denote the corresponding "decomposition" by
\[\pi = \pi_{p} = \setsuchas{s \cup\set{0}}{s \in S}.\]

The simplest case is when \(\pi=\set{\NN}\).
The resulting convolution is the Dirichlet convolution, with \(A(n)=\divisors{n}\).
Cashwell and Everett \cite{NumThe} famously proved that \(\Gamma\) with this product
is isomorphic to the (full) ring of power series on countably many indeterminates
(for the \(i\)'th prime number \(p_{i}\), the indicator function \(e_{{p_{i}}}\) correspond to the indeterminate \(x_{i}\)), and
that this ring is a unique factorization domain. This is in line with
Remark \ref{orgb245134}, since in this case only the prime powers \(p^{1}\) are primitive.
The dense subring generated by the \(e_{p}\)'s is the ring of finitely supported arithmetical functions.
This subring is isomorphic to the ring of polynomials in countably many indeterminates.

We will henceforth restrict our study to \emph{bounded} homogeneous regular convolutions, by which
we mean convolutions where there exists an \(M < \infty\)
such that all "blocks" \(u \in \pi\) have block-length \(\le M\). The Dirichlet convolution is
obviously not bounded.

\begin{definition}
By a slight abuse of notation, for a bounded homogeneous regular convolution we say that
\begin{enumerate}[label={(\alph*)}]
\item the type of \(a \in \NNp\) is the type of \(p^{a}\), where \(p\) is any prime number,
\item \(a\) is primitive whenever \(p^{a}\) is,
\item for a primitive \(a\),  the rank of \(a\) is the rank of \(p^{a}\),
\item if \(n=ka\) with \(a\) primitive and \(1 \le k \le d\), we say that \(n\) has \emph{height} \(k\)  and write \(k=\height{d}{n}\).
\end{enumerate}
\label{org70b1361}
\end{definition}

\begin{remark}
If \(M \subset \NNp\) denotes the primitive elements (w.r.t. some bounded homogeneous regular convolution), then
\begin{equation}
\label{eq:orgb39f0db}
\NNp = \bigsqcup_{{k=1}}^{d} kM,
\end{equation}
where the union is disjoint, and
\[
kM = \setsuchas{km}{m \text{ primitive}} = \setsuchas{n \in \NNp}{\height{d}{n}=k}.\]

Thus, if \(p\) is a prime with corresponding block
\[
\pi_{p} = \set{0, b, 2b, 3b, \dots, sb}
\]
then the non-zero elements of the block have type \(b\),
\(b\) is the only primitive element,
the rank of \(b\) is \(s\),
and the height of \(kb\) is \(k\).
As mentioned in the introduction, we say that \(\pi_{p}\) has block-length \(s\).
\end{remark}
\subsection{The unitary convolution}
\label{sec:org1cfd398}
The simplest example of a bounded homogeneous regular convolution is when all blocks \(u \in \pi\) have
block-length 1. Then \(S\) is the partition
of \(\NNp\) into singleton sets. All prime powers are primitive, of rank one. This yields the
unitary convolution (\ref{eq:org182c43e}).
Note that if \(p,q\) are different primes, and \(a,b \in \NNp\), then
\begin{equation}
\label{eq:orge654700}
  \begin{split}
    e_{p^{a}} * e_{q^{b}} & = e_{p^a q^b,} \\
    e_{p^{a}} * e_{p^{b}} & = \mathbf{0}.    
  \end{split}
\end{equation}
\section{Ternary convolution}
\label{sec:orgeaefe31}
\subsection{Definition via set-partition}
\label{sec:orga2923a5}
The ternary convolution (\ref{eq:org9b5b66e}) has already been introduced,
but we will now define it in another way.

\begin{theorem}
There is a unique regular homogeneous bounded convolution where all block-lengths
are equal to 2. This \emph{ternary convolution} is defined by a single partition
of \(\NNp\) into parts \(\set{m,2m}\), which is as follows. Put 
\begin{equation}
\label{eq:org002fe6e}
M =\setsuchas{4^{k} \ell}{\ell \text{ odd}} = \setsuchas{m \in \NNp}{\val_{2}(n) \text{ even}}.
\end{equation}
Then
\[ 2M = \setsuchas{2m}{m \in M} = \setsuchas{n \in \NNp}{\val_{2}(n) \text{ odd}}, \]

and
\[\NNp = M \, \bigsqcup \, 2M,\] where the union is disjoint.
The parts of the partition of \(\NNp\) associated with the convolution are
\begin{equation}
\label{eq:org7573984}
\setsuchas{\set{m,2m}}{m \in M}
\end{equation}
(and the "blocks" are \(\set{0,m,2m}\)). Thus for \(p\) prime and \(m \in M\), \(p^{m}\) is primitive, of
rank two, type \(m\), and \(p^{2m}\) is not primitive, of type \(m\).
\label{org8fb3da9}
\end{theorem}
\begin{proof}
If we choose to \textbf{not} place a number of the form \(v=2^{s}d\) with \(s,d\) odd
in a block together with \(v/2\), then \(v/2\) ends up in a block \(\set{0,v/2}\) of
block-length 1.
\end{proof}

The first parts of the partition has already been shown in 
Table \ref{tab:orge53dce4}.
A few more are displayed in
Table \ref{tab:org297cd36}.
Reading the integers "block by block" we obtain \seqnum{A036552}.

\begin{table}[ptbh]
\caption{\label{tab:org297cd36}Parts of the partition defining the ternary convolution.}
\centering
\begin{tabular}{l|l|l|l|l|l|l|l|l|l}
part \# & 1 & 2 & 3 & 4 & 5 & 6 & 7 & 8 & 9\\
\hline
part & 1, 2 & 3, 6 & 4, 8 & 5, 10 & 7, 14 & 9, 18 & 11, 22 & 12, 24 & 13, 26\\
\end{tabular}
\end{table}
\subsection{Multiplication table of the primitives}
\label{sec:org2664034}
The multiplicative structure of \((\Gamma,*,+)\), with \(*\) the ternary convolution,
is determined from the multiplication table of \(e_{{p^{m}}}\)'s with \(p^{m}\) primitive.
Here we have that
\begin{equation}
\label{eq:orgc7842dd}
\begin{split}
  e_{p^{m}} * e_{p^{m}} &= e_{p^{2m}} \\
  e_{p^{m}} * e_{p^{m}} * e_{p^{m}} = e_{p^{m}}^{3} & =  \mathbf{0} 
\end{split}
\end{equation}
If \(p,q\) are distinct primes, and \(a,b\) positive integers, then
\(e_{{p^{a}}} * e_{{q^{b}}} = e_{{p^{a}}{q^{b}}}\).
\subsection{Natural density}
\label{sec:org13dc707}
Recall the notion of asymptotic, or natural, density (see for instance the "Density" section of the wiki \cite{oeis}):
\begin{definition}
If \(S \subset \NN\) then the \emph{lower asymptotic density}  of \(S\) is
\[\liminf_{{n \to \infty}} \frac{\# \setsuchas{s \in S}{s \le n}}{n},\]
and the \emph{upper asymptotic density}  of \(S\) is
\[\limsup_{{n \to \infty}} \frac{\# \setsuchas{s \in S}{s \le n}}{n}.\]
If these limits coincide, then the \emph{natural density} (or \emph{asymptotic density}) \(\dens(S)\) exists and is
equal to this common value.
\end{definition}

The following  properties of the natural density are elementary:
\begin{lemma}
Let \(S,T \subset \NNp\) both have natural density, and let \(k \in \NNp\).
Then  \(\dens(S + \set{k}) = \dens(S)\), and 
\(\dens(S) = k\dens(kS)\). If \(S\) is finite, then \(\dens(S)=0\), and
\(\dens(\NNp \setminus S)=1\). More generally, if \(S \cap T = \emptyset\), then \(\dens(S \cup T) =\dens(S) + \dens(T)\).
\label{org560ea21}
\end{lemma}

\begin{lemma}
Let \((s_{i})_{i=0}^{\infty}\) be a sequence of positive integers such that \(s_{i+1} -s_{i} > 1\) for all \(i\).
Let \(p\) be a prime number, let \(S\) denote the image of the sequence, and put
\[
T = \setsuchas{n \in \NNp}{\val_{p}(n) \in S}.
\]
Then \(T\) has natural density
\begin{equation}
\label{eq:org0201447}
\dens(T) = \sum_{k=0}^{\infty}(p^{-s_{k}} - p^{-s_{k}-1}).
\end{equation}
\label{org5fa90b4}
\end{lemma}
\begin{proof}
It is well-known \cite{Conroy:Density} that the set
\[\setsuchas{n \in \NNp}{\val_{p}(n) = s_{k}}\]
has natural density
\(p^{-s_{k}} - p^{-s_{k}-1}.\)
Furthermore, it follows immediately from Lemma \ref{org560ea21} that if
\(A_{1},\dots,A_{r}\) are finitely many disjoint sets, having natural density \(c_{1},\dots,c_{r}\),
then \(\cup_{i=1}^{r} A_{i}\) has natural density \(\sum_{{i=1}}^{r}c_{i}\).
Since the right-hand side of (\ref{eq:org0201447}) converges, the result follows.
\end{proof}

\begin{corollary}
Let \(p\) be a prime number, \(n>1\)  an integer, \(d \ge 0\) an integer.
Then the set
\[
\setsuchas{m \in \NNp}{\val_p(m) \equiv \modd{d}{n}}
\]
has natural density
\[
  \sum_{k=0}^{\infty} (p^{-nk -d} - p^{-nk -d -1}) =
  p^{-d}(1-p^{-1})\sum_{k=0}^{\infty} p^{-nk} =
  p^{-d}(1-p^{-1}) (1-p^{-n})^{-1}.
\]
\label{orgefa6660}
\end{corollary}

Similarly, we have that
\begin{corollary}
Let \(p \neq q\) be  prime numbers, \(n,m>1\)  integers, \(c,d \ge 0\) integers.
Then the set
\[
\setsuchas{s \in \NNp}{\val_p(s) \equiv \modd{d}{n}, \, \val_q(s) \equiv \modd{c}{m}}
\]
has natural density
\[
  p^{-d}(1-p^{-1}) (1-p^{-n})^{-1}    q^{-c}(1-q^{-1}) (1-q^{-m})^{-1.}
\]
\label{orgf3d0d82}
\end{corollary}
\begin{proof}
Follows from Corollary \ref{orgefa6660} and the Chinese Remainder Theorem.
\end{proof}
\subsection{Density of primitives for the ternary convolution}
\label{sec:org6f7d53b}
We now consider the the set \(M\) defined by (\ref{eq:org002fe6e}), in other words, the primitive exponents for the ternary convolution.
It starts as
\begin{table}[ptbh]
\caption{\label{tab:orgbc05f78}Primitives of the ternary convolutions.}
\centering
\begin{tabular}{rrrrrrrrrrrrrrrrrrrr}
1 & 3 & 4 & 5 & 7 & 9 & 11 & 12 & 13 & 15 & 16 & 17 & 19 & 20 & 21 & 23 & 25 & 27 & 28 & 29\\
\end{tabular}
\end{table}

This is sequence \seqnum{A003159} in OEIS \cite{oeis}, the so-called "vile numbers". Conroy
\cite{Conroy:Density} uses this very sequence to illustrate the concept of natural density.
Lemma \ref{orgefa6660} shows that \(M\) has density \(2/3\). Of course, if we knew in advance
that \(M\) has a natural density \(c\), we could use the fact that \(\NNp = M \sqcup 2M\) to
conclude that \(1 =c + c/2\), and so \(c=2/3\).

\begin{figure}[ptbh]
\centering
\includegraphics[width=0.7 \textwidth]{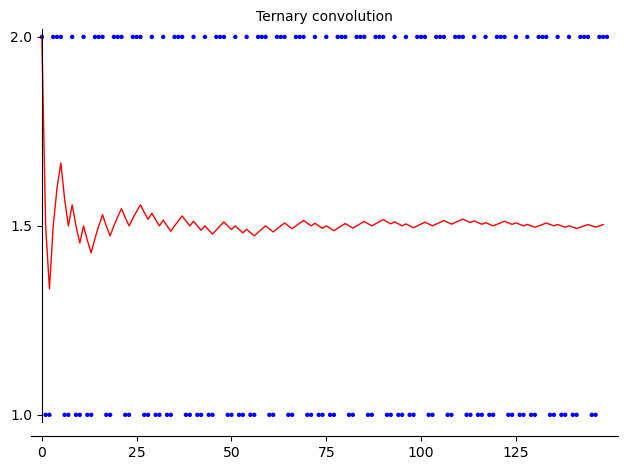}
\caption{\label{fig:orgda2a76e}Distances (in blue) between consecutive primitives for the ternary convolution; cumulative average in red. Horizontal axis indexed by ordinality, so \(n\) means \(n\)'th primitive exponent.}
\end{figure}

The distance between two consecutive primitive elements is either one or two, since any odd number is primitive.
Furthermore, an even number cannot be followed by an even number.
Figure \ref{fig:orgda2a76e} plots this distance, and the cumulative averages.
Since the natural density of the primitives is \(2/3\),
the \(n\)'th primitive is approximately \(\frac{3}{2}n\)  and the
\((n+1)\)'th primitive is approximately \(\frac{3}{2}(n+1)\), so the average distance is approximately \(3/2\).

\section{Greedy convolutions}
\label{sec:org603cab4}
\subsection{Definition of greedy convolutions: trees and matrices}
\label{sec:orgd1cec38}
\begin{definition}
Let \(d\) be a positive integer. Define the family of rooted trees \(T_{d}(n)\), \(n \in \NN\), where the vertices of \(T_{d}(n)\) are the integers \(\set{0,1,\dots, n}\),
by the following:
\begin{enumerate}
\item Let \(T_{d}(0) = 0\), i.e., the rooted tree with a single vertex 0.
\item Let \(T_{d}(n+1)\) be obtained from \(T_{d}(n)\) by adding the vertex \(n+1\) as a leaf to a branch
\[0 \to k \to 2k \to \cdots \to mk,\]
extending it to
\[0 \to k \to 2k \to \cdots \to mk \to (m+1)k,\]
where
\[n+1= (m+1)k, \quad m+1 \le d.\]
\item If there are several such branches, choose the one with the smallest \(k\).
\item If there is no such branch, attach \(n+1\) to the root \(0\), starting a new branch.
\end{enumerate}

Then:
\begin{enumerate}
\item Define \(T_{d}\) as the inductive limit (union) of the \(T_{d}(n)\)'s.
\item Let \(\gdecomp{d}\) be the decomposition of \(\NN\) into arithmetic progressions containing 0 that one
obtains by considering the branches of \(T\) as the blocks.
\item Let \(\gpartition{d}\) be the corresponding
partition of \(\NNp\), obtained by removing the root 0 from \(T\) and considering the components of the
so obtained forest as the parts of \(\gpartition{d}\).
\item Define \(\gconvolution{d}\), the \textbf{greedy convolution} of (block-)length \(d\), as the homogeneous regular convolution given by the decomposition
\(\gdecomp{d}\) (or equivalently by the partition \(\gpartition{d}\)).
\end{enumerate}
\label{org58460e1}
\end{definition}

\begin{figure}[ptbh]
\centering
\includegraphics[width=0.5 \textwidth]{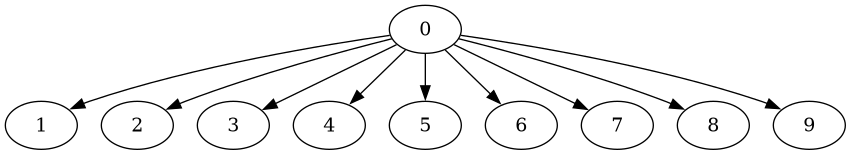}
\caption{\label{fig:org6587207}The tree \(T_{1}(9)\).}
\end{figure}

\begin{figure}[ptbh]
\centering
\includegraphics[width=0.5 \textwidth]{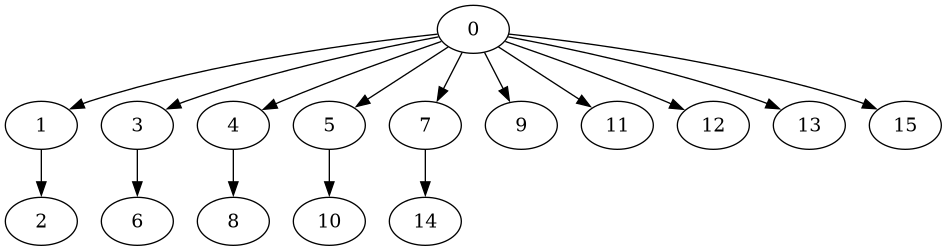}
\caption{\label{fig:org45d0380}The tree \(T_{2}(15)\).}
\end{figure}

\begin{figure}[ptbh]
\centering
\includegraphics[width=0.9 \textwidth]{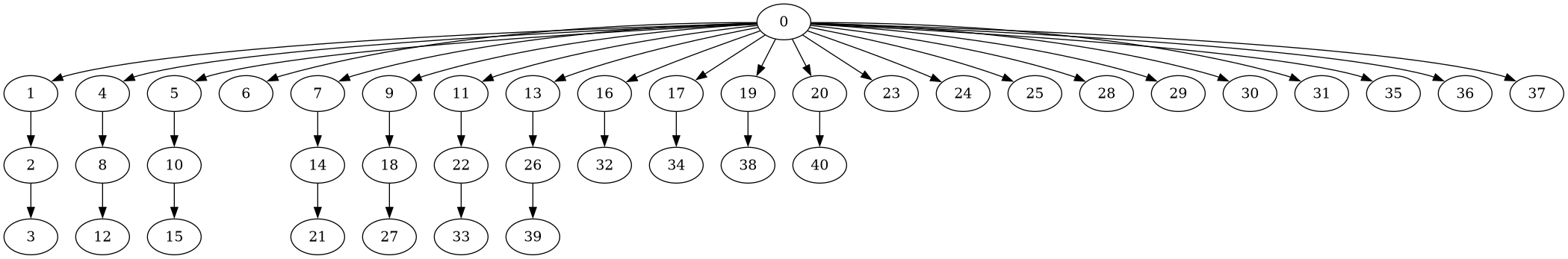}
\caption{\label{fig:orga3a2aae}The tree \(T_{3}(40)\).}
\end{figure}

\begin{figure}[ptbh]
\centering
\includegraphics[width=0.9 \textwidth]{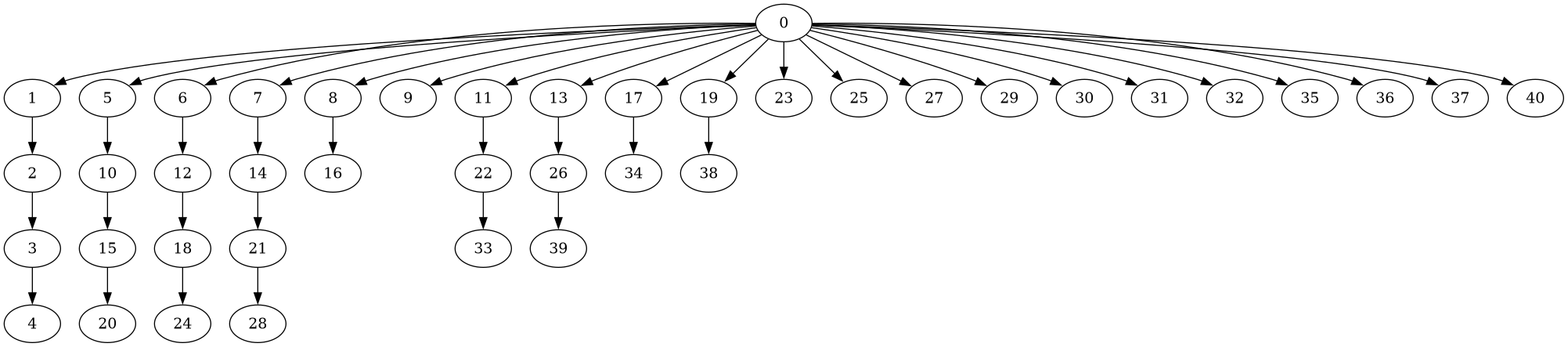}
\caption{\label{fig:org0f37182}The tree \(T_{4}(40)\).}
\end{figure}

We show \(T_{1}(9)\), \(T_{2}(15)\), \(T_{3}(40)\), and \(T_{4}(40)\)
in Figures \ref{fig:org6587207} - \ref{fig:org0f37182}.

\begin{lemma}
The partition \(\gdecomp{d}\) can be obtained as follows:
\begin{enumerate}
\item Form the \(d \times N\)-matrix \(A^{(d,N)}=(a_{i,j})\) with \((a_{i,j}) = ij\).
\item Starting with column \(c=2\), loop through the elements \(a_{r,c}\) with \(r\) from 1 to \(d\).
\item If the value \(a_{r,c}\) occurs in \(A^{(d,N)}\) in some earlier column, put \(a_{\ell,c}=0\) for all \(\ell \ge r\).
\item Continue until the last column.
\item Optionally, delete the all-zero columns.
\item The non-zero columns of \(A^{(d,N)}\) now contains those parts of \(\gpartition{d}\) which have primitive element
(i.e., smallest element in the part) \(\le N\). The primitive elements are the non-zero elements in the first row.
\end{enumerate}
\label{orga5162d2}
\end{lemma}
\begin{proof}
Both descriptions capture the same idea: the integer \(n\) can be assigned to the part with smallest element \(k\)
if \(n=k \ell\) with \(1 \le \ell \le d\). If it fits into different parts, put it into the part with smallest \(k\).
\end{proof}

\begin{mexample}
Consider the tree \(T_{3}(40)\) shown in Figure \ref{fig:orga3a2aae}. The branch \(0 \to 6\) will never be added to in any \(T_{4}(n)\),
since 12 is already taken. Thus 6 will be in a part of its own in \(\gpartition{3}\).
\end{mexample}

\begin{lemma}
The greedy convolutions of length 1 and 2 are  the unitary convolution and the ternary convolution.
\label{org97675d4}
\end{lemma}
\begin{proof}
The unitary convolution is the unique homogeneous regular convolution where all blocks in \(\pi\) have size two,
and the ternary convolution is the unique homogeneous regular convolution where all blocks have size three.
\end{proof}

\begin{remark}
The unitary convolution similarly has \(e_{p^{m}}^{2}=0\) for all prime powers,
and parts of length one. It is the greedy convolution of length 1.
The "ternary convolution" was named "ternary" for the following reasons: \(e_{p^{m}}^{3}=0\) for all prime powers,
the blocks of the decomposition \(\pi_{p}\) have size 3, (but block-length 2),
and the parts of the corresponding partition have cardinality 2.
The ternary convolution is the greedy convolution of block-length 2. Maybe the "dyadic convolution" would
be a better name, but we will stick to the name that Gavell (and the present author, who
was the Bachelor's thesis advisor for \cite{Gavell:Regular}) coined.
\end{remark}

\begin{mexample}
The first greedy convolutions have  the following matrices, before pruning the all-zero columns.
To start, the unitary convolution has singleton blocks:

\begin{displaymath}
A^{(1,20)} =
\left(\begin{array}{rrrrrrrrrrrrrrrrrrrr}
1 & 2 & 3 & 4 & 5 & 6 & 7 & 8 & 9 & 10 & 11 & 12 & 13 & 14 & 15 & 16 & 17 & 18 & 19 & 20
\end{array}\right)
\end{displaymath}

The ternary convolution has no unfilled positions (all non-zero columns have no zeroes):
\begin{displaymath}
A^{(2,20)} =
\left(\begin{array}{rrrrrrrrrrrrrrrrrrrr}
1 & 0 & 3 & 4 & 5 & 0 & 7 & 0 & 9 & 0 & 11 & 12 & 13 & 0 & 15 & 16 & 17 & 0 & 19 & 20 \\
2 & 0 & 6 & 8 & 10 & 0 & 14 & 0 & 18 & 0 & 22 & 24 & 26 & 0 & 30 & 32 & 34 & 0 & 38 & 40
\end{array}\right)
\end{displaymath}

The greedy convolution of length three has non-zero columns with one or three non-zero elements.
In other words, the primitives have rank one, or full rank 3.
\begin{displaymath}
A^{(3,20)} =
\left(\begin{array}{rrrrrrrrrrrrrrrrrrrr}
1 & 0 & 0 & 4 & 5 & 6 & 7 & 0 & 9 & 0 & 11 & 0 & 13 & 0 & 0 & 16 & 17 & 0 & 19 & 20 \\
2 & 0 & 0 & 8 & 10 & 0 & 14 & 0 & 18 & 0 & 22 & 0 & 26 & 0 & 0 & 32 & 34 & 0 & 38 & 40 \\
3 & 0 & 0 & 12 & 15 & 0 & 21 & 0 & 27 & 0 & 33 & 0 & 39 & 0 & 0 & 48 & 51 & 0 & 57 & 60
\end{array}\right)
\end{displaymath}

The greedy convolution of length 4 has primitives of rank one, rank 2, and full rank 4, but seemingly no
primitives of rank 3.
\begin{displaymath}
A^{(4,20)} =
\left(\begin{array}{rrrrrrrrrrrrrrrrrrrr}
1 & 0 & 0 & 0 & 5 & 6 & 7 & 8 & 9 & 0 & 11 & 0 & 13 & 0 & 0 & 0 & 17 & 0 & 19 & 0 \\
2 & 0 & 0 & 0 & 10 & 12 & 14 & 16 & 0 & 0 & 22 & 0 & 26 & 0 & 0 & 0 & 34 & 0 & 38 & 0 \\
3 & 0 & 0 & 0 & 15 & 18 & 21 & 0 & 0 & 0 & 33 & 0 & 39 & 0 & 0 & 0 & 51 & 0 & 57 & 0 \\
4 & 0 & 0 & 0 & 20 & 24 & 28 & 0 & 0 & 0 & 44 & 0 & 52 & 0 & 0 & 0 & 68 & 0 & 76 & 0
\end{array}\right)
\end{displaymath}

The greedy convolution of length 5 is more complicated still.
\begin{displaymath}
A^{(5,20)} =
\left(\begin{array}{rrrrrrrrrrrrrrrrrrrr}
1 & 0 & 0 & 0 & 0 & 6 & 7 & 8 & 9 & 10 & 11 & 0 & 13 & 0 & 15 & 0 & 17 & 0 & 19 & 0 \\
2 & 0 & 0 & 0 & 0 & 12 & 14 & 16 & 0 & 20 & 22 & 0 & 26 & 0 & 0 & 0 & 34 & 0 & 38 & 0 \\
3 & 0 & 0 & 0 & 0 & 18 & 21 & 0 & 0 & 0 & 33 & 0 & 39 & 0 & 0 & 0 & 51 & 0 & 57 & 0 \\
4 & 0 & 0 & 0 & 0 & 24 & 28 & 0 & 0 & 0 & 44 & 0 & 52 & 0 & 0 & 0 & 68 & 0 & 76 & 0 \\
5 & 0 & 0 & 0 & 0 & 30 & 35 & 0 & 0 & 0 & 55 & 0 & 65 & 0 & 0 & 0 & 85 & 0 & 95 & 0
\end{array}\right)
\end{displaymath}
\end{mexample}
\subsection{Simple properties}
\label{sec:org3b7ff7e}

\begin{lemma}
Consider the greedy convolution of length \(d\), and let \(n\) be a positive integer.
\begin{enumerate}[label={(\roman*)}]
\item The type of \(n\) is the element on the same branch of \(T_{d}(n)\) which is adjacent to the root. Equivalently, it is the
element in the same column but in row one of \(A^{(d)}\).
\item The integer \(n\) is primitive if and only if it is adjacent to the root in \(T_{d}(n)\), in other words,
if it is in row one of \(A^{(d)}\).
\item If \(n\) is primitive, its rank is the length of the branch of \(T_{d}\) that contains it; equivalently,
it is the number of non-zero elements in the column of \(A^{(d,N)}\) that contains it (for \(N\) sufficiently large).
\item The height \(\height{d}{n}\) of \(n\) is the distance in \(T_{d}(n)\) from  the root to \(n\) (equivalently,
the index of the row in \(A^{(d,n)}\) that contains \(n\)).
\end{enumerate}
\label{org5efffdc}
\end{lemma}
\begin{proof}
Immediate; note that the length of a branch of \(T_{d}(n)\) is the number of edges in it, or equivalently the
distance from the leaf to the root.
\end{proof}

\begin{lemma}
Consider the greedy convolution of length \(d\).
\begin{enumerate}[label={(\roman*)}]
\item The exponent \(1\) is primitive, and \(2,3,\dots,d\) are not primitive.
\item If \(n\) is not primitive, there is a unique pair
\[(a,k) \in \set{1,2,\dots,n-1} \times \set{2,\dots,d}\] with
\(n=ak\) and \(a\) primitive; then \(\height{d}{n}=k\).
\item The exponent \(n\) is \textbf{not primitive} if and only if
\(n=ak\) with \(a\) primitive and \(1 < k \le d\), and
\[\height{d}{2a}=2, \height{d}{3a}=3, \dots, \height{d}{(k-1)a}=k-1.\]
\end{enumerate}
\label{org7ffedaf}
\end{lemma}
\begin{proof}
We can assume that \(n > d\). The first branch of \(T_{d}(n)\) is always
\[0 \to 1 \to 2 \to \cdots \to d,\]
showing
that \(1\) is primitive and \(2,3,\dots,d\) are not primitive.
The number \(n\) is a leaf of \(T_{d}(n)\), sitting on some branch
\[0 \to a \to 2a \to \cdots \to ka=n,\]
where
\(a\) is primitive. It is always the case that \(k \le d\). The integer \(n\) has been placed at its "spot"
because there were no earlier available "spots", where a "spot" is an empty position \(X\) at
\[0 \to b \to 2b \to \cdots \to (c-1)b \to X,\]

so that \(n=cb\) with \(c \le d\), \(b<a\) primitive. If \(k=1\) then
\(n\) is primitive, otherwise it is not, and has \(\height{d}{n}=k\).
\end{proof}

\begin{remark}
It is \textbf{not the case} that the integer \(n>d\) is  primitive if and only if
there is no \(1 < k \le d\) such that \(\divides{k}{n}\) and  \(\frac{n}{k}\) is primitive.
As an example, for \(d=4\) the integer 32 is primitive even though \(32/4=8\) is primitive.
\end{remark}
\subsection{Uniqueness of the unitary and the ternary convolution}
\label{sec:org8d3cb47}
\begin{theorem}[Gavell]
For \(d > 2\), the partition \(\gpartition{d}\) contains some part \(s \in \gpartition{d}\) with
\(\#s < d\).
\label{org1191dce}
\end{theorem}
\begin{proof}
Suppose not, and that the first few blocks are
\begin{equation}
\label{eq:org688900a}
\begin{split}
b_{1} & = \set{0, 1, 2, \dots, d} \\
b_{2} & = \set{0, d+1, 2(d+1), 3(d+1), \dots, d(d+1)} \\
b_{3} & = \set{0, d+2, 2(d+2), 3(d+2), \dots, d(d+2)} \\
b_{4} & = \set{0, d+3, 2(d+3), 3(d+3), \dots, d(d+3)} \\
b_{5} & = \set{0, d+4, 2(d+4), 3(d+4), \dots, d(d+4)} \\
   &  \qquad \vdots \\
  b_{k} & = \set{0, s, 2s, \dots, ds}
\end{split}
\end{equation}
where \(s\) is the smallest number not contained in a previous parts.

If \(d\) is odd, then consider the blocks \(b_{2}\) and \(b_{4}\) starting with \(d+1\) and \(d+3\) as their
first non-zero element. Since
\[\gcd(d+1,d+3)=2\] we get that
\[w=\lcm(d+1,d+3)=\frac{(d+1)(d+3)}{2}.\]
Since \(d \ge 3\), the inequality \(\frac{d+3}{2} < d-1\) follows; hence \(w\) could fit into both
\(b_{2}\) and \(b_{4}\).
The greedy procedure places it in \(b_{2}\), so block \(b_{4}\) will never fill up:
it consists of only the \(\frac{d}{2}\) elements
\[\set{0,d+2, 2(d+2), \dots ,(\frac{d}{2})(d+3)}.\]
As shown in Figure \ref{fig:orga3a2aae}, for \(d=3\) the
positive integer 12 has already been added to the block \(b_{2}=\set{0,4,8,12}\).
Therefore, the block \(b_{4} = \set{0,6}\) cannot be filled.

If \(d\) is even, then consider the blocks \(b_{3}\) and \(b_{5}\) starting with \(d+2\) and \(d+4\), respectively,
as their smallest positive element. We now have that \[\gcd(d+2,d+4)=2\] and that
\[w=\lcm(d+2,d+4) = \frac{(d+2)(d+4)}{2}.\]
Since \(d \ge 4\) we have that \(\frac{d+4}{2} < d\), so once again
\(w\) fits into two different blocks, \(b_{3}\) and \(b_{5}\).
In Figure \ref{fig:org0f37182} we illustrate that for \(d=4\), the integer 24 can fit into both
the block starting with \(6\) and the block starting with \(8\), leaving the latter block unfilled.
\end{proof}

\begin{remark}
The proof is taken almost verbatim from \cite{Gavell:Regular},
but we replaced \(d-1\) with \(d\), since we choose to index
the convolutions by the size of the \emph{parts} (or the length of the blocks) rather than the size of the blocks.

A more careful analysis shows that when \(d\) is odd, the block with smallest non-zero element \(d+3\)
is in fact the \textbf{first} unfilled block; when \(d\) is even it is the block starting with \(d+4\).

The proof rules out the existence of a regular homogeneous convolution where all parts of of the
partition of \(\NNp\) have equal size \(d \ge 3\),
since if you choose not to be greedy, you immediately get an unfilled block.
\label{orgd2651e8}
\end{remark}

We conclude:
\begin{corollary}
The unitary convolution and the ternary convolution are the only homogeneous bounded regular convolutions where all blocks have the same length.
\label{org6ef57c1}
\end{corollary}
\section{The greedy convolution of block-length 3}
\label{sec:org1f75020}

\subsection{Primitive elements}
\label{sec:org81de2f9}
We tabulate the parts of the greedy convolution of block-length 3 in Table \ref{tab:org063fdbe}.
The primitive elements, first those of full rank and then those of rank one, are listed in Table \ref{tab:org3771568}.
It turns out that these primitive elements are related to sequence
 \seqnum{A339690}.

\begin{table}[ptbh]
\caption{\label{tab:org063fdbe}Parts of \(\gpartition{3}\)}
\centering
\begin{tabular}{lllll}
(1, 2, 3) & (4, 8, 12) & (5, 10, 15) & (6,) & (7, 14, 21)\\
(9, 18, 27) & (11, 22, 33) & (13, 26, 39) & (16, 32, 48) & (17, 34, 51)\\
(19, 38, 57) & (20, 40, 60) & (23, 46, 69) & (24,) & (25, 50, 75)\\
(28, 56, 84) & (29, 58, 87) & (30,) & (31, 62, 93) & (35, 70, 105)\\
(36, 72, 108) & (37, 74, 111) & (41, 82, 123) & (42,) & (43, 86, 129)\\
(44, 88, 132) & (45, 90, 135) & (47, 94, 141) & (49, 98, 147) & (52, 104, 156)\\
(53, 106, 159) & (54,) & (55, 110, 165) & (59, 118, 177) & (61, 122, 183)\\
(63, 126, 189) & (64, 128, 192) & (65, 130, 195) & (66,) & (67, 134, 201)\\
(68, 136, 204) & (71, 142, 213) & (73, 146, 219) & (76, 152, 228) & (77, 154, 231)\\
(78,) & (79, 158, 237) & (80, 160, 240) & (81, 162, 243) & (83, 166, 249)\\
(85, 170, 255) & (89, 178, 267) & (91, 182, 273) & (92, 184, 276) & (95, 190, 285)\\
(96,) & (97, 194, 291) & (99, 198, 297) & (100, 200, 300) & (102,)\\
\end{tabular}
\end{table}

\begin{table}[ptbh]
\caption{\label{tab:org3771568}Primitive elements of rank 3 and of rank one of \(\gpartition{3}\)}
\centering
\begin{tabular}{rrrrrrrrrrrrrrr}
1 & 4 & 5 & 7 & 9 & 11 & 13 & 16 & 17 & 19 & 20 & 23 & 25 & 28 & 29\\
31 & 35 & 36 & 37 & 41 & 43 & 44 & 45 & 47 & 49 & 52 & 53 & 55 & 59 & 61\\
63 & 64 & 65 & 67 & 68 & 71 & 73 & 76 & 77 & 79 & 80 & 81 & 83 & 85 & 89\\
91 & 92 & 95 & 97 & 99 & 100 & 101 & 103 & 107 & 109 & 112 & 113 & 115 & 116 & 117\\
- & - & - & - & - & - & - & - & - & - & - & - & - & - & -\\
6 & 24 & 30 & 42 & 54 & 66 & 78 & 96 & 102 & 114 & 120 & 138 & 150 & 168 & 174\\
186 & 210 & 216 & 222 & 246 & 258 & 264 & 270 & 282 & 294 & 312 & 318 & 330 & 354 & 366\\
\end{tabular}
\end{table}

\begin{theorem}
For the partition \(\gpartition{3}\) associated with the greedy convolution of length 3, we have that:
\begin{enumerate}[label={(\roman*)}]
\item The part sizes of \(\gpartition{3}\) are 3 or 1.
\item The initial elements in the parts of size 3 are the positive integers
\[M= \setsuchas{m \in \NNp}{m=4^{i}9^{j}k, \, \gcd(k,6)=1}.\]
\item For every \(m \in M\) there is a singleton part \(\set{6m}\). There are no other singleton parts.
\item The numbers \(M\) constitute OEIS \seqnum{A339690} and the numbers \(3M\)
constitute OEIS \seqnum{A329575}.
\item The natural density of \(M\) is \(1/2\),
thus the natural density of \(6M\) is \(1/12\), and
the natural density of
\[\setsuchas{a \in \NNp}{p^{a} \text{ is primitive }}\]
 is \(7/12\).
\end{enumerate}
\end{theorem}
\begin{proof}
The comments to OEIS \seqnum{A339690} states that \(M\) is the set of positive integers obtained by successively adding
\(1,2,3,\dots\), but only adding \(n\) if it is not \(2a\), \(3a\), or \(6a\) for some \(a\) already added.
Such a number \(n\) starts a new branch of \(T_{3}(n)\), which will be filled to
\[0 \to m \to 2m \to 3m.\]
Then \(6n\) will start a new branch, but this branch will never grow further, since
\(12n\) will attach to
\[0 \to 4n \to 8n \to 12n.\]
It follows that level one of \(T_{4}\) consists of \(M\) and \(6M\),
level two consists of \(2M\), and level three of \(3M\).
Hence the positive integers are a disjoint union
\[\NNp = M \,\sqcup \, 2M \, \sqcup \, 3M \, \sqcup \, 6M,\]
and Corollary \ref{orgf3d0d82} shows that \(M\) has natural density.
Let \(c\) be the natural density of \(M\).  Then
\[1 = c + \frac{c}{2} + \frac{c}{3} + \frac{c}{6},\]
which implies that \(c=\frac{1}{2}\).

Hence \(6M\) has natural density \(\frac{1}{12}\) and the primitives,
which are \(M \sqcup 6M\), have natural density \(\frac{1}{2} + \frac{1}{12} = \frac{7}{12}\).
\end{proof}
\subsection{Statistics of the greedy convolution of length 3}
\label{sec:orgce26e4f}
\begin{figure}[ptbh]
\centering
\includegraphics[width=0.7 \textwidth]{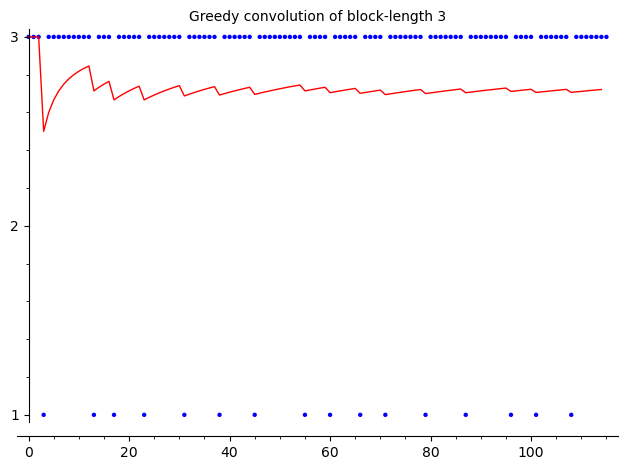}
\caption{\label{fig:orga2146a4}Length of parts of greedy convolution of length 3 (in blue), cumulative average in red. Horizontal axis indexed by ordinality, so \(n\) means \(n\)'th primitive exponent.}
\end{figure}

In figure \ref{fig:orga2146a4} we plot the lengths of the parts in \(\gpartition{3}\).
From our density results we have that the ratio of parts of size three to parts of size one is six to one.
Thus the average part-size should be
\[\frac{1}{7} \times 1 + \frac{6}{7} \times 3 = \frac{19}{7} \approx 2.7,\]
which is hinted at by the red line in the graph, showing cumulative averages of block-lengths of the first blocks.

\begin{figure}[ptbh]
\centering
\includegraphics[width=0.9 \textwidth]{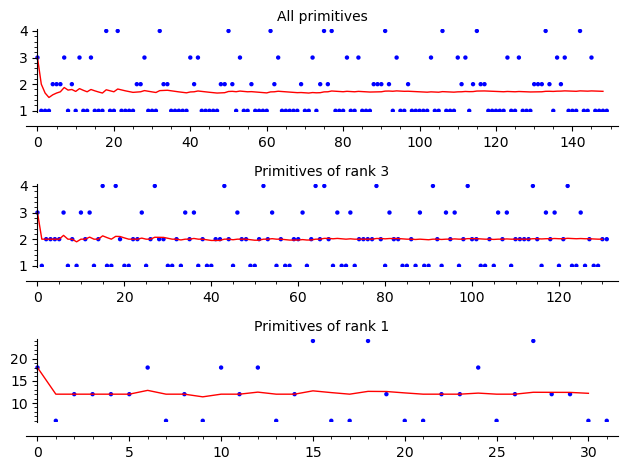}
\caption{\label{fig:orga394d10}Distances (in blue) between consecutive primitives for greedy convolution of length 3; cumulative averages in red. Horizontal axis indexed by ordinality, so \(n\) means \(n\)'th primitive exponent, or \(n\)'th primitive exponent of rank 1, or \(n\)'th primitive exponent of rank 3.}
\end{figure}

The density of the primitives being 7/12, it follows that the average distance between
to primitives is \(12/7 \approx 1.71\). The primitives of order \(3\) have density \(1/2\), so the average distance between
them is 2; finally the primitives of order \(1\) have density \(1/12\), so the average distance between
them is 12. This is shown in figure \ref{fig:orga394d10}.
\section{Selective sifting}
\label{sec:org6b6b851}
Using the method of \emph{selective sifting}, described below, 
we can describe the primitive elements for the greedy convolution with block-length 2,
and the primitive elements of full rank for the greedy convolution of block-length 3.

\begin{definition}
Let \(S \subset \NNp \setminus \set{1}\). Define the \emph{selective sifting} \(\ssift{S} \subset \NNp\) inductively by
\begin{equation}
  \begin{split}
    \ssift{S}_{\le 1} & = \set{1} \\
    \ssift{S}_{\le n+1} & =
\begin{cases}
\ssift{S}_{\le n},  & n+1 \in \setsuchas{ab}{a \in S, \, b \in \ssift{S}_{\le n}}; \\
\ssift{S}_{\le n} \cup \set{n+1}, & \text{otherwise.}
\end{cases} \\
    \ssift{S} & = \cup_{n=1}^{\infty} \ssift{S}_{\le n} 
  \end{split}
\end{equation}
\label{org0285b2c}
\end{definition}

\begin{lemma}
\(T=\ssift{S}\) is the unique subset of \(\NNp\) such that
for \(n \in \NNp \setminus \set{1}\),
\[
  n \in T \quad \iff \quad n \not \in \setsuchas{ab}{a \in S, \, b \in T,  b < n}.
\]
\end{lemma}

\begin{remark}
Informally, we sift out multiples of elements of \(S\), but only multiples by factors that have note already been sifted out.
In particular, if \(n\) is not divisible by any \(s \in S\) then \(n \in \ssift{S}\), but if \(n \in S\) then \(n \not\in \ssift{S}\).
\end{remark}

\begin{mexample}[The sifting of $S=\{2\}$]
\begin{table}[tbhp]
\caption{\label{tab:orgd454ebf}Values of \(\mathbf{s}(\{2\})_{\le k}\)}
\centering
\begin{tabular}{rl}
k & s(\{2\})\textsubscript{<= k}\\
1 & \{1\}\\
2 & \{1\}\\
3 & \{1, 3\}\\
4 & \{1, 3, 4\}\\
5 & \{1, 3, 4, 5\}\\
6 & \{1, 3, 4, 5\}\\
7 & \{1, 3, 4, 5, 7\}\\
8 & \{1, 3, 4, 5, 7\}\\
9 & \{1, 3, 4, 5, 7, 9\}\\
10 & \{1, 3, 4, 5, 7, 9\}\\
11 & \{1, 3, 4, 5, 7, 9, 11\}\\
12 & \{1, 3, 4, 5, 7, 9, 11, 12\}\\
13 & \{1, 3, 4, 5, 7, 9, 11, 12, 13\}\\
14 & \{1, 3, 4, 5, 7, 9, 11, 12, 13\}\\
15 & \{1, 3, 4, 5, 7, 9, 11, 12, 13, 15\}\\
16 & \{1, 3, 4, 5, 7, 9, 11, 12, 13, 15, 16\}\\
17 & \{1, 3, 4, 5, 7, 9, 11, 12, 13, 15, 16, 17\}\\
18 & \{1, 3, 4, 5, 7, 9, 11, 12, 13, 15, 16, 17\}\\
\end{tabular}
\end{table}

In Table \ref{tab:orgd454ebf} we compute the first steps of the sifting for \(S=\set{2}\).
After infinitely many steps, we obtain \(\ssift{\set{2}}\). Combining
Theorem \ref{org8fb3da9} with the next Lemma, we get that \(\ssift{\set{2}}\)
is nothing but the primitives for the ternary convolution,
shown in Table \ref{tab:orgbc05f78}.
\label{orge7d6215}
\end{mexample}

\begin{lemma}
If \(p\) is prime, \(S=\set{p}\), \(M=\ssift{S}\) then
\[M = \setsuchas{n \in \NNp}{\val_{p}(n) \text{ even}},\]
and
\[\NNp = M \sqcup pM.\]
The natural density of \(M\) is
\(\frac{1}{1+ \frac{1}{p}}\).
\label{orgdad3ece}
\end{lemma}
\begin{proof}
We use induction on \(n\), the induction base being \(n=1 \in\ssift{p}\). Next assume that \(n>1\).
If \(n=p^{{2k}}r\) with \(\gcd(p,r)=1\) then
\[n/p = p^{{2k-1}}r \not \in \ssift{p},\]
by the induction hypothesis, so
\(n \in \ssift{p}\).
If \(n=p^{{2k+1}}r\) with \(\gcd(p,r)=1\) then
\[n/p = p^{{2k}}r \in \ssift{p},\]
so \(n \not \in \ssift{p}\).

To show that \(\NNp = M \sqcup pM\), note that if \(n=p^{a}r\) with \(\gcd(p,r)=1\) then \(pn = p^{a+1}r\); furthermore,
\(a\) is either even or odd.

Corollary \ref{orgefa6660} gives that \(M\) has natural density and that this density is
\(\frac{1}{1+ \frac{1}{p}}\).
\end{proof}

\begin{remark}
The sequences \(\ssift{\set{5}}\) and \(\ssift{\set{7}}\) are
\seqnum{A382744} and \seqnum{A382745}.
\end{remark}

\begin{lemma}
If \(p,q\) are distinct primes, \(S=\set{p,q,pq}\), \(M=\ssift{S}\) then
\[M =\setsuchas{n \in \NNp}{\val_{p}(n) \text{ even and } \val_{q}(n) \text{ even}},\]
and
\[\NNp = M \sqcup pM \sqcup qM \sqcup pqM.\]
The natural density of \(M\) is
\(\frac{1}{1 + \frac{1}{p} + \frac{1}{q} + \frac{1}{pq}}\).
\label{orge402740}
\end{lemma}
\begin{proof}
We first prove that
\(M =\setsuchas{n \in \NNp}{\val_{p}(n) \text{ even and } \val_{q}(n) \text{ even}}\).
We use induction on \(n\), the base case being \(n=1\).
For \(n>1\), write \(n=p^{a}q^{b}r\) with \(\gcd(pq,r)=1\).
We want to prove that \(n \in \ssift{S}\) if and only if \(a,b\) are both even.
\begin{itemize}[label={\(\star\)}]
\item If \(a\) is even, \(b\) is odd then \(n/q = p^{a}q^{b-1}r \in \ssift{S}\) by the induction hypothesis, so \(n \not \in \ssift{S}\).
\item If \(a\) is odd, \(b\) is even then \(n/p = p^{a-1}q^{b}r \in \ssift{S}\), so \(n \not \in \ssift{S}\).
\item If \(a\) is odd, \(b\) is odd then \(n/(pq) = p^{a-1}q^{b-1}r \in \ssift{S}\), so \(n \not \in \ssift{S}\).
\item If \(a=b=0\) then no element of \(S\) divides \(n\), so \(n \in \ssift{S}\).
\item If \(a=0\), \(b>0\) is even then \(n/q = q^{b-1}r \not \in \ssift{S}\), so \(n \in \ssift{S}\).
\item If \(b=0\), \(a>0\) is even then \(n/p = p^{a-1}r \not \in \ssift{S}\), so \(n \in \ssift{S}\).
\item If \(a>0\) is even, \(b>0\) is even then
\begin{align*}
n/p & = p^{a-1}q^{b}r \not \in \ssift{S},\\
n/q & = p^{a}q^{b-1}r \not \in \ssift{S}, \\
n/(pq) &  = p^{a-1}q^{b-1}r \not \in \ssift{S},
\end{align*}
so \(n \in \ssift{S}\).
\end{itemize}

We have  proved that \(n \in \ssift{S}\) if and only if \(a,b\) are both even.  

Next, we show that  \(\NNp = M \sqcup pM \sqcup qM \sqcup pqM\).
Note that if  \(n=p^{a}q^{b}r\) with \(\gcd(pq,r)=1\) then
\begin{align*}
pn &= p^{a+1}q^{b}r, \\
qn &= p^{a}q^{b+1}r, \\
pqn &= p^{a+1}q^{b+1}r,
\end{align*}
which covers all combinations of \(a,b\) modulo 2.

The assertion about the natural density of \(M\) follows from Corollary \ref{orgf3d0d82}.
Alternatively, since Corollary \ref{orgf3d0d82} guarantees
that the natural density \(c\) exists,
\[1= c + \frac{c}{p} + \frac{c}{q} + \frac{c}{pq},\]
so
\[
c=\frac{1}{1 + \frac{1}{p} + \frac{1}{q} + \frac{1}{pq}}.
\]
\end{proof}

\begin{mexample}
The first elements of \(\ssift{\set{2,3,6}}\) are
\begin{displaymath}
{1, 4, 5, 7, 9, 11, 13, 16, 17, 19, 20, 23, 25, 28, 29, 31, 35, 36, 37}
\end{displaymath}
and Lemma \ref{orge402740} shows that \(\ssift{\set{2,3,6}}\) consists of those positive integers
with an even number of two's and an even number of three's in their prime factorization;
this coincides with the primary elements of rank 3 for the greedy convolution of length 3.
\end{mexample}

We give some more "selective sifting" results that are not needed to describe the greedy convolutions
we have studied so far, but are amusing in their own right.

\begin{lemma}
If \(p\) is a prime, \(k\) a positive integer, \(S=\set{p^{k}}\), \(M=\ssift{S}\) then
\[M = \setsuchas{n \in \NNp}{\val_{p}(n) \equiv 0,1,\dots,k-1 \hspace{-2mm} \pmod{2k}} .\]
\label{org4000ed9}
\end{lemma}
\begin{proof}
The assertion is true for \(n=1\). We prove it by induction.
For  \(n>1\), write
\[n=p^{a}r, \qquad \gcd(p,r)=1.\]
If \(a < k\) then \(\dividesnot{p^{k}}{n}\), so \(n \in S\).

If \(a \ge k\) then write
\(a=2ck + d,  0 \le d < 2k\).

There are then two subcases. If
\[a=2ck + d, \qquad 0 \le d < k,\] then
\[n/p^{k} = p^{a'}r, \qquad a' = 2(c-1)k + d+k,\]
so by the induction hypothesis \(n/p^{k} \not \in M,\) hence \(n \in M\).

If
\[a=2ck + d, \qquad k \le d < 2k,\]  then
\[n/p^{k} = p^{a'}r, \qquad a' = 2ck + d-k,\]
so by the induction hypothesis \(n/p^{k}  \in M,\) hence \(n \not \in M\).
\end{proof}

\begin{remark}
The sequence \(\ssift{\set{8}}\) is \seqnum{A382746} in OEIS.
\end{remark}

\begin{theorem}
If \(p\) is a prime, \(k\) a positive integer, \(S=\set{p,p^{k}}\), \(M=\ssift{S}\) then if \(k\) is odd
\[M = \setsuchas{n \in \NNp}{\val_{p}(n) \equiv \modd{0}{2}},\]
and if \(k\) is even then
\[M = \setsuchas{n \in \NNp}{\val_{p}(n) \equiv \modd{0,2,4, \dots, k-2}{k+1}}.\]
\label{org2dba548}
\end{theorem}
\begin{proof}
Assume first that \(k\) is odd.
We want to prove that
\[\ssift{\set{p,p^{k}}} = \setsuchas{n \in \NNp}{\val_{p}(n) \equiv \modd{0}{2}}.\]

Write \(n=p^{a}r\) with \(\gcd(p,r)=1\). We proceed by induction on \(k\).

If \(a=0\) then neither \(p\) nor \(p^{k}\) divides \(n\), so \(n \in \ssift{\set{p,p^{k}}}\).

If \(a>0\) then \(n/p=p^{a-1}r\).
If \(a\) is even then \(a-1\) is odd,
so \(n/p \not \in \ssift{\set{p,p^{k}}}\), by the induction hypothesis.
If \(a\) is odd then \(a-1\) is even,
so \(n/p \in \ssift{\set{p,p^{k}}}\).

If \(a \ge k\) then \(n/p^{k}=p^{a-k}r\).
If \(a\) is even then \(a-k\) is odd,
so \(n/p^{k} \not \in \ssift{\set{p,p^{k}}}\), by the induction hypothesis.
If \(a\) is odd then \(a-k\) is even,
so \(n/p^{k} \in \ssift{\set{p,p^{k}}}\).

Thus for all \(d \in \set{p,p^{k}}\) we see that if \(\divides{d}{n}\) then
\(n/d \in \ssift{\set{p,p^{k}}}\) if and only if \(a\) is odd. Hence
\(n=p^{a}r \in \ssift{\set{p,p^{k}}}\) if and only if \(a\) is even, as stated.

Next assume that \(k>0\) is even. Write
\[n=p^{a}r, \qquad \gcd(p,r)=1,\]
and write
\[a= c(k+1) + d,  \qquad 0 \le d < k+1.\]
Then \(n/p = p^{a-1}r\) (if \(a>0\)) and \(n/p^{k} = p^{a-k}r\) (if \(a \ge k\)).
Write the set of congruence classes modulo \(k+1\) as a disjoint union
\begin{displaymath}
\frac{\ZZ}{(k+1)\ZZ}  = A \sqcup B
\end{displaymath}
with
\begin{align*}
A & =  \set{[0]_{k+1}, [2]_{k+1}, [4]_{k+1}, \dots, [k-2]_{{k+1}}} \\
B &=  \set{[1]_{k+1}, [3]_{k+1}, [5]_{k+1}, \dots, [k-1]_{k+1}, [k]_{k+1}}
\end{align*}
We want to prove that
\[n=p^{a}r \in \ssift{p,p^{k}} \qquad \iff \qquad [a]_{k+1} \in A.\]

Observe that
\begin{itemize}[label={\(\star\)}]
\item Subtracting \([1]_{k+1}\) from an element of \(A\) yields an element of \(B\).
\item Subtracting \([k]_{{k+1}}\) from an element of \(A\) yields an element of \(B\).
\item Subtracting \([k]_{k+1}\) from an element of \(B\) yields an element of \(A\).
\item Subtracting \([1]_{k+1}\) from an element of \(B \setminus \set{[k]_{k+1}}\) yields an element of \(A\);
subtracting \([k]_{k+1}\) from \([k]_{k+1}\) yields \([0]_{k+1} \in A\).
\end{itemize}

With a little book-keeping we see that if
\[n=p^{a}r, \qquad [a]_{k+1} \in A\]
then if \(1 \le a\) then
\[n/p = p^{a-1}r, \qquad  [a-1]_{k+1} \in B,\]
hence
\[n/p \not \in \ssift{\set{p,p^{k}}}.\]
Similarly,  if
\(a \ge k\) then
\[n/p^{k} = p^{a-k}r, \qquad [a-k]_{k+1} \in B,\]
hence
\[n/p^{k} \not \in \ssift{\set{p,p^{k}}}.\]
In both cases, there is no allowed factor in \(\ssift{p,p^{k}}\), so
\(n \in \ssift{p,p^{k}}\).

Conversely,
if \[n=p^{a}r,\qquad [a]_{k+1} \in B,\] then
if \(1 \le a\) with \(a \not \in [k]_{k+1}\) then
\[n/p = p^{a-1}r,\qquad [a-1]_{k+1} \in A,\]
hence
\[n/p  \in \ssift{p,p^{k}} \quad  \text{ so } \quad n  \not \in \ssift{p,p^{k}}.\]
Finally, if \(a \ge k\) (this includes the case \(a  \in [k]_{k+1}\))
then \[n/p^{k} = p^{a-k}r,\qquad [a-k]_{k+1} \in A,\] hence
\[n/p^{k} \in \ssift{p,p^{k}} \quad \text{ so } \quad n  \not \in \ssift{p,p^{k}}.\]
\end{proof}

\begin{theorem}
If \(p,q\) are distinct primes, \(S=\set{p,q}\), \(M=\ssift{S}\) then
\[M =\setsuchas{n \in \NNp}{\val_{p}(n) \equiv \modd{\val_{q}(n)}{2}}.\]
\end{theorem}
\begin{proof}
Induction, similar to the proof of Theorem \ref{org2dba548}.
\end{proof}
\subsection{Questions and conjectures regarding selective sifting}
\label{sec:org791a854}
\begin{question}
Let \(n\) be an integer \(>1\).
What is \(\ssift{\set{n}}\), and what is its natural density?
\end{question}

\begin{question}
Let \(p\) be prime, and let \(a_{1},\dots,a_{k}\) be positive integers.
What is \[\ssift{\set{p^{a_{1}}, \dots,p^{a_{1}}}}\] and what is its natural density?
\end{question}

\begin{question}
If \(p,q\) are distinct primes, \(S_{1} = \set{p}\), \(S_{2} = \set{q}\), \(S_{3} = \set{p,q,pq}\),
then natural densities of \(\ssift{\set{S_{i}}}\) are \(c_{1}=\frac{1}{1+1/p}\),
\(c_{2}=\frac{1}{1+1/q}\), \(c_{3}=c_{1}c_{2}\). How can this be generalized?
\end{question}

\begin{question}
What are the possible natural densities achievable by \(\ssift{\set{S}}\) when \(S\) has
\(k\) elements?
\end{question}
\section{The greedy convolution of length 4}
\label{sec:orgfd29dc1}
The non-zero columns of \(A^{(4)}\) starts as follows:

\begin{displaymath}
\left(\begin{array}{rrrrrrrrrrrrrrr}
1 & 5 & 6 & 7 & 8 & 9 & 11 & 13 & 17 & 19 & 23 & 25 & 27 & 29 & 30 \\
2 & 10 & 12 & 14 & 16 & 0 & 22 & 26 & 34 & 38 & 46 & 50 & 54 & 58 & 60 \\
3 & 15 & 18 & 21 & 0 & 0 & 33 & 39 & 51 & 57 & 69 & 75 & 81 & 87 & 90 \\
4 & 20 & 24 & 28 & 0 & 0 & 44 & 52 & 68 & 76 & 92 & 100 & 108 & 116 & 120
\end{array}\right)
\end{displaymath}

This table is \seqnum{A382747}.

There are primitives of rank one, two and four. The primitives of rank 4 are listed in Table \ref{tab:orgc9a0059},
together with the ones of rank one and two. The concatenation
of these lists, i.e., the sequence of all primitives, is \seqnum{A382748}.

\begin{table}[ptbh]
\caption{\label{tab:orgc9a0059}Primitives of rank four, one, and two for the greedy convolution of length 4}
\centering
\begin{tabular}{rrrrrrrrrrrrrrr}
1 & 5 & 6 & 7 & 11 & 13 & 17 & 19 & 23 & 25 & 27 & 29 & 30 & 31 & 32\\
35 & 37 & 41 & 42 & 43 & 47 & 49 & 53 & 55 & 59 & 61 & 65 & 66 & 67 & 71\\
73 & 77 & 78 & 79 & 83 & 85 & 89 & 91 & 95 & 97 & 101 & 102 & 103 & 107 & 109\\
113 & 114 & 115 & 119 & 121 & 125 & 127 & 131 & 133 & 135 & 137 & 138 & 139 & 143 & 144\\
145 & 149 & 150 & 151 & 155 & 157 & 160 & 161 & 162 & 163 & 167 & 169 & 173 & 174 & 175\\
179 & 181 & 185 & 186 & 187 & 189 & 191 & 193 & 197 & 199 & 203 & 205 & 209 & 210 & 211\\
215 & 217 & 221 & 222 & 223 & 224 & 227 & 229 & 233 & 235 & 239 & 241 & 245 & 246 & 247\\
251 & 253 & 256 & 257 & 258 & 259 & 263 & 265 & 269 & 271 & 275 & 277 & 281 & 282 & 283\\
287 & 289 & 293 & 294 & 295 & 297 & 299 & 301 & 305 & 307 & 311 & 313 & 317 & 318 & 319\\
323 & 325 & 329 & 330 & 331 & 335 & 337 & 341 & 343 & 347 & 349 & 351 & 352 & 353 & 354\\
355 & 359 & 361 & 365 & 366 & 367 & 371 & 373 & 377 & 379 & 383 & 385 & 389 & 390 & 391\\
395 & 397 & 401 & 402 & 403 & 407 & 409 & 413 & 415 & 416 & 419 & 421 & 425 & 426 & 427\\
431 & 433 & 437 & 438 & 439 & 443 & 445 & 449 & 451 & 455 & 457 & 459 & 461 & 462 & 463\\
467 & 469 & 473 & 474 & 475 & 479 & 481 & 485 & 487 & 491 & 493 & 497 & 498 & 499 & 503\\
505 & 509 & 510 & 511 & 513 & 515 & 517 & 521 & 523 & 527 & 529 & 533 & 534 & 535 & 539\\
541 & 544 & 545 & 546 & 547 & 551 & 553 & 557 & 559 & 563 & 565 & 569 & 570 & 571 & 575\\
- & - & - & - & - & - & - & - & - & - & - & - & - & - & -\\
9 & 45 & 48 & 63 & 99 & 117 & 153 & 171 & 207 & 216 & 225 & 240 & 243 & 261 & 279\\
- & - & - & - & - & - & - & - & - & - & - & - & - & - & -\\
8 & 36 & 40 & 56 & 88 & 104 & 136 & 152 & 180 & 184 & 192 & 200 & 232 & 248 & 252\\
280 & 296 & 328 & 344 & 376 & 392 & 396 & 424 & 440 & 468 & 472 & 488 & 520 & 536 & 568\\
\end{tabular}
\end{table}

\begin{conjecture}
Let \(M\) denote the primitives of rank 4 for the greedy convolution of length 4.
\begin{itemize}
\item The primitives of rank one and two are all divisible by 8 or 9.
\item "Most" primitives of rank one or two belong to \(8M\) or \(9M\).
\item There are no primitives of rank 3.
\end{itemize}
\end{conjecture}
\begin{remark}
Note that, for instance \(36=9 \times 4\), yet \(4\) is not a primitive of rank 4.
\end{remark}

If we try to express the \(M\),
the set of primary elements of rank 4 for the greedy convolution of length 4,
using selective sifting, we should successively add elements to \(S\) that belong to
\(\ssift{S} \setminus M\), and hope that \(M \setminus \ssift{S}\) "takes care of itself".
However, when sifting best as possible the integers \(\le N=4000\), we get that
\begin{equation}
\label{eq:org7f2dcc3}
\begin{split}
S & = \set{2,3,4,8,9,16,36,72,192,384,768,1024,2048} \\
(M \cap[N]) \setminus (\ssift{S} \cap [N]) & = \emptyset \\
(\ssift{S} \cap [N]) \setminus (M \cap[N])  & = \set{256,1280,1792,2816,3328}
\end{split}
\end{equation}
This means that \(M\) \textbf{is not} expressible using selective sifting. However,
it seems "close" to the selective sifting \(\ssift{S}\) with an infinite \(S\)
that has some structure:
\begin{conjecture}
Let \(M\) denote the primitives of rank 4 for the greedy convolution of length 4,
and let \(S\) be the "greedy sifting set" where we start by \(S=\set{2,3,4}\) and
successively augment it by adding the smallest element of
\(M \setminus \ssift{S}\).
Then \(S\) will become infinite, and its elements will all be of the form \(2^{a}3^{b}\).

The "sporadic" elements in
\(\ssift{S}  \setminus M\) will all be of the form \(2^{8}p\) where \(p\) is an odd prime; this set
is either finite or very sparse.
\end{conjecture}
\section{Questions about greedy convolutions}
\label{sec:org44d0849}

\subsection{Possible ranks of primitive elements}
\label{sec:orgf2cbb70}
\begin{lemma}
If \(d > 2\) is odd, then the first primitive number of rank one in the greedy convolution of block-length \(d\)
is \(3\frac{d+1}{2}\); if \(d\) is even, the first such is \(3 \frac{d+2}{2}\).
\end{lemma}
\begin{proof}
The number \(1\) is primitive and \(2,3,\dots,d\) are not. Then \(d+1,\dots,2d\) are primitive, by Lemma (\ref{org7ffedaf}).
In the odd case, one part is
\[\set{d+1, 2(d+1),3(d+1),\dots}\] and another is
\(\set{3\frac{d+1}{2}}\), and cannot contain \(2 \times 3\frac{d+1}{2}\), since it is already taken.

In the even case, one part is \[\set{d+2, 2(d+2),3(d+2),\dots}\] and another is
\(\set{3\frac{d+2}{2}}\), and cannot contain \(2 \times 3\frac{d+2}{2}\), since it is already taken.
\end{proof}

We tabulate the first primitives of a given rank for the greedy convolutions of block-length  \(\le 16\) in
Table \ref{tab:orgfb58e13}. A zero indicates the we have not found a primitive of the given type.

\begin{question}
Consider the greedy convolution of block-length \(d>2\).
Are there primitives of rank \(d-1\)?
\end{question}

\begin{question}
More generally, what are the ranks of primitives?
\end{question}

\begin{question}
What are the smallest primitives of a rank type?
\end{question}

\begin{question}
What is the proportion and distribution of the ranks of the primitives?
\end{question}

\begin{table}[ptbh]
\caption{\label{tab:orgfb58e13}First primitive of given rank for greedy convolutions of block-length 2 to 16}
\centering
\begin{tabular}{rrrrrrrrrrrrrrrrr}
Len & Rk &  &  &  &  &  &  &  &  &  &  &  &  &  &  & \\
 & 1 & 2 & 3 & 4 & 5 & 6 & 7 & 8 & 9 & 10 & 11 & 12 & 13 & 14 & 15 & 16\\
2 & 0 & 1 &  &  &  &  &  &  &  &  &  &  &  &  &  & \\
3 & 6 & 0 & 1 &  &  &  &  &  &  &  &  &  &  &  &  & \\
4 & 9 & 8 & 0 & 1 &  &  &  &  &  &  &  &  &  &  &  & \\
5 & 9 & 8 & 40 & 0 & 1 &  &  &  &  &  &  &  &  &  &  & \\
6 & 12 & 80 & 10 & 72 & 0 & 1 &  &  &  &  &  &  &  &  &  & \\
7 & 12 & 80 & 10 & 35 & 504 & 0 & 1 &  &  &  &  &  &  &  &  & \\
8 & 15 & 12 & 105 & 14 & 0 & 96 & 0 & 1 &  &  &  &  &  &  &  & \\
9 & 15 & 16 & 105 & 12 & 896 & 729 & 0 & 0 & 1 &  &  &  &  &  &  & \\
10 & 18 & 16 & 15 & 162 & 14 & 63 & 567 & 0 & 0 & 1 &  &  &  &  &  & \\
11 & 18 & 16 & 15 & 162 & 14 & 63 & 567 & 616 & 99 & 0 & 1 &  &  &  &  & \\
12 & 21 & 20 & 125 & 18 & 0 & 16 & 144 & 99 & 0 & 0 & 0 & 1 &  &  &  & \\
13 & 21 & 20 & 125 & 18 & 0 & 16 & 144 & 99 & 0 & 143 & 0 & 0 & 1 &  &  & \\
14 & 24 & 20 & 125 & 18 & 0 & 352 & 22 & 99 & 0 & 143 & 0 & 0 & 0 & 1 &  & \\
15 & 24 & 28 & 20 & 168 & 21 & 180 & 18 & 99 & 0 & 143 & 0 & 0 & 0 & 0 & 1 & \\
16 & 27 & 24 & 25 & 32 & 21 & 112 & 0 & 20 & 880 & 143 & 0 & 195 & 0 & 0 & 0 & 1\\
\end{tabular}
\end{table}
\subsection{Natural density of primitives}
\label{sec:org900c641}
\begin{figure}[ptbh]
\centering
\includegraphics[width=0.9\textwidth]{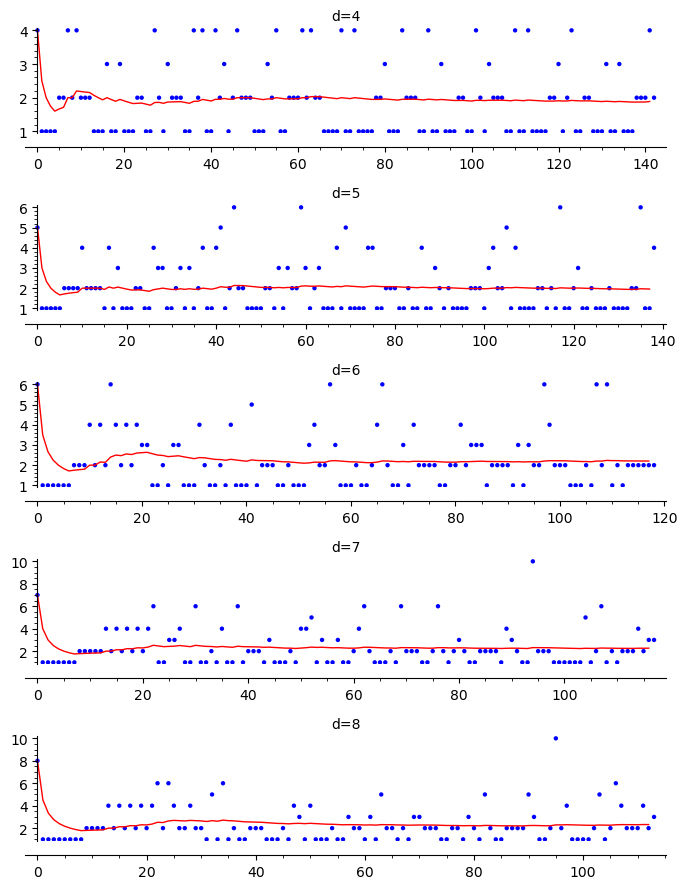}
\caption{\label{fig:orgffe7a12}Distance between consecutive primitives for greedy convolution of length \(4,5,6,7,8\). Cumulative averages in red. Horizontal axis indexed by ordinality, so \(n\) means \(n\)'th primitive exponent.}
\end{figure}

As noted for the greedy convolution of block-length 3, the average distance between two consecutive primitives
should be the reciprocal of the natural density of the primitives. We plot these distances and their cumulative averages
in Figure \ref{fig:orgffe7a12}.

\begin{question}
What is the natural density of the primitives for the greedy convolution of block-length \(d\)?
\end{question}

\subsection{Height of elements}
\label{sec:orgce6072b}
\begin{figure}[ptbh]
\centering
\includegraphics[width=0.9\textwidth]{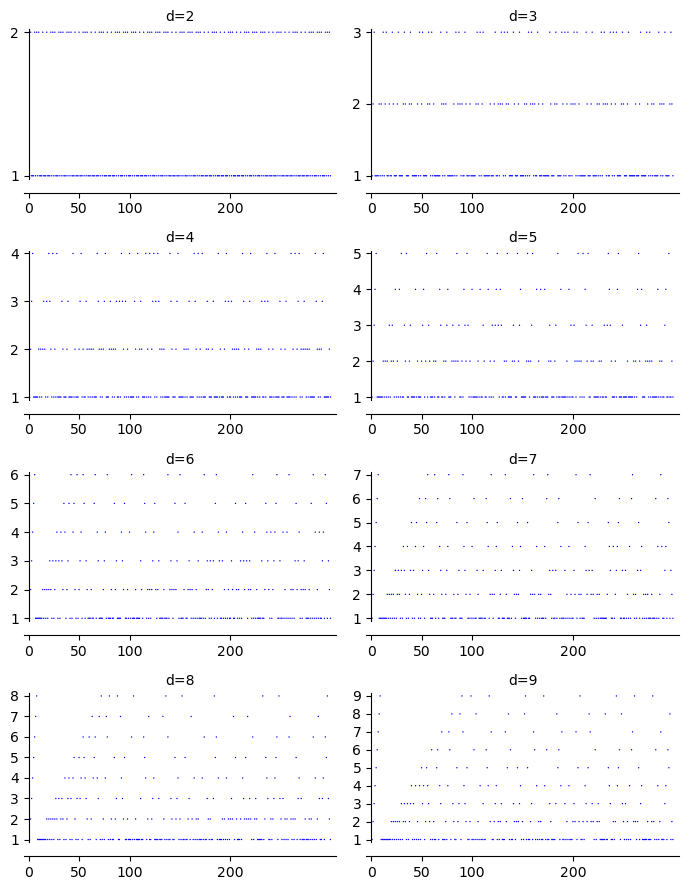}
\caption{\label{fig:org9178ce8}Heights of elements with respect to greedy convolutions of length \(2,3,4,5,6,7,8,9\). Horizontal axis indexed naturally, so \(n\) means \(n\).}
\end{figure}

The heights of the first 300 positive integers with respect to
the greedy convolutions of block-length \(2,3,4,5,6,7,8,9\) are as plotted in Figure \ref{fig:org9178ce8}.

The heights of integers with respect to the greedy convolution of block-length 4 is \seqnum{A382749}.

\begin{question}
Apart from the obvious fact that \(h_{d}(n) = n\) for \(1 \le n \le d\), there seems to be sporadic lines of other slopes
appearing. How can these features be explained?
\end{question}
\begin{question}
For the greedy convolution of block-length \(d\), what is the density of integers of height \(k\)?
\end{question}

\section{AI disclosure}
The first version of this manuscript, \href{https://arxiv.org/abs/2504.02795}{arXiv:2504.02795}, deposited on April 2025,
was written by the author with no AI usage (but relying heavily on the SageMath
\cite{sagemath} computer algebra system). The same is true for all revisions excluding this one (deposited on August 14, 2026). The impetus for this version was that an adversarial review by Claude (Anthropic) discovered errors
in Definition 8 and Lemma 15
(in both cases, \(<\) should be corrected to \(\le\))
and in the proof of Theorem 14
(\(\frac{c}{12}\) corrected to \(\frac{c}{6}\))
and in Corollary 7 (\(p\) corrected to \(q\)).

Edits to the layout of the corrected manuscript, and in particular the \emph{layout} of the bibliography (not its contents)
were likewise made by the LLM agent.
It also compiled ancillary items, such as
\begin{itemize}
  \item  a short document listing and describing the  OEIS sequences occurring in the manuscript,
\item   a document explaining the use of the attached SageMath code.      
\end{itemize}
Note that the code, included as \detokenize{anc/ternary_code.sage},
is written solely by the author, without LLM assistance (which partly explains the code's clunkiness).

No new results or material were provided by the LLM agent, and the text of the manuscript is written by the author without AI aid, except for some \LaTeX\ code tweaks, as explained above.

\bibliographystyle{jis}
\bibliography{greedy}

\end{document}